\documentclass{baustms}
\citesort

\theoremstyle{cupthm}
\newtheorem{thm}{Theorem}[section]

\newtheorem{lemma}[thm]{Lemma}
\theoremstyle{cupdefn}

\theoremstyle{cuprem}
\newtheorem{rem}[thm]{Remark}
\numberwithin{equation}{section}
\newtheorem{conj}[thm]{Conjecture}

\begin{document}
\runningtitle{Euclid-Euler Heuristics}
\title{Euclid-Euler Heuristics for Perfect Numbers}
%% If there is more than one author, put \cauthor immediately before
%% the corresponding author.
%\cauthor %% mark the next author as corresponding author
\author[1]{Jose Arnaldo B. Dris}
\address[1]{PhD Student, UP - Diliman, Quezon City, PH\email{josearnaldobdris@gmail.com, jadris@feu.edu.ph}}
%% If there are several authors, list them here
%\author[2]{Second author}
%\address[2]{Second address\email{a@net.com}}

%% List the authors, initials and surnames only, for the
%% running head (left hand page)
\authorheadline{J.~A.~B.~Dris}

%% If there is a dedication, include it here
%\dedication{Dedicated to ...}

\dedication{Dedicated to Dr. Severino V. Gervacio, for his suggestion of tackling Sorli's conjecture first}

\support{The author does not currently receive support from any academic institution.}

\begin{abstract}
An odd perfect number $N$ is said to be given in \emph{Eulerian form} if $N = {q^k}{n^2}$ where $q$ is prime with $q \equiv k \equiv 1 \pmod 4$ and $\gcd(q,n) = 1$.  Similarly, an even perfect number $M$ is said to be given in \emph{Euclidean form} if $M = (2^p - 1)\cdot{2^{p - 1}}$ where $p$ and $2^p - 1$ are primes.  In this article, we show how simple considerations surrounding the differences between the underlying properties of the Eulerian and Euclidean forms of perfect numbers give rise to what we will call the \emph{Euclid-Euler heuristics} for perfect numbers.
\end{abstract}

%% - subject classification and keywords
%% 2010 American Mathematical Society Subject Classification
%% Provide only ONE primary classification
%% Four or five keywords or phrases

\maketitle

\section{Introduction}
If $J$ is a positive integer, then we write $\sigma(J)$ for the sum of the divisors of $J$.  We denote the abundancy index of $x \in \mathbb{N}$ as $I(x)=\sigma(x)/x$, where $\mathbb{N}$ is the set of natural numbers or positive integers.  A number $L$ is \emph{perfect} if $\sigma(L)=2L$.

An even perfect number $M$ is said to be given in \emph{Euclidean form} if $$M = (2^p - 1)\cdot{2^{p - 1}}$$
where $p$ and $2^p - 1$ are primes.  We call $M_p = 2^p - 1$ the \emph{Mersenne prime} factor of $M$.  Currently, there are only $49$ known Mersenne primes, which correspond to
$49$ even perfect numbers.

An odd perfect number $N$ is said to be given in Eulerian form if $$N = {q^k}{n^2}$$
where $q$ is prime with $q \equiv k \equiv 1 \pmod 4$ and $\gcd(q,n) = 1$.  We call $q^k$ the \emph{Euler part} of $N$ while $n^2$ is called the \emph{non-Euler part} of $N$.

It is currently unknown whether there are infinitely many even perfect numbers, or whether any odd perfect numbers exist.  It is widely believed that there is an infinite number of even perfect numbers. On the other hand, no examples for an odd perfect number have been found (despite extensive computer searches), nor has a proof for their nonexistence been established.

Ochem and Rao \cite{OchemRao} recently proved that $N > {10}^{1500}$.  Nielsen \cite{Nielsen1} has obtained the lower bound $\omega(N) \geq 10$ for the number of \emph{distinct} prime factors of $N$, improving on his last result $\omega(N) \geq 9$ (see \cite{Nielsen2}).

Sorli conjectured in \cite{Sorli} that $k=\nu_{q}(N)=1$ always holds.  Dris conjectured in \cite{Dris3} and \cite{Dris5} that the divisors $q^k$ and $n$ are related by the inequality $q^k < n$.  This conjecture was predicted on the basis of the result $I(q^k)<\sqrt[3]{2}<I(n)$.

Broughan, et.~al.~ \cite{BroughanDelbourgoZhou} recently showed that for any odd perfect number $N = {q^k}{n^2}$, the ratio of the non-Euler part $n^2$ to the Euler part $q^k$ is greater than $315/2$.  This improves on a result of Dris \cite{Dris3}.

In a recent paper, Chen and Chen \cite{ChenChen} improves on Broughan, et.~al.'s results, and poses a related (open) problem.

\section{The General Multiplicative Form of All Perfect Numbers}\label{Section1}
Suppose that $N = {q^k}{n^2}$ is an odd perfect number given in Eulerian form.  Since prime powers are deficient and $\gcd(q, n) = 1$, we know that $q^k \neq n$.  (In particular, it is also true that $q \neq n$.)  Consequently, we know that either $q^k < n$ or $n < q^k$ is true.

Observe that the Euclidean form $M = (2^p - 1)\cdot{2^{p - 1}}$ for an even perfect number $M$ possesses a multiplicative structure that is \emph{almost} similar to that of the Eulerian form $N = {q^k}{n^2}$ for an odd perfect number $N$.  Here is a table comparing and contrasting the underlying properties of these two forms, which we shall refer to as the \emph{Euclid-Euler heuristics} for (odd) perfect numbers:

\begin{center}
\begin{tabular}{|l|l|}
\hline
(E-1) (Euclid-Euler Theorem)  & (O-1) (Conjecture, 2010 \cite{Dris4}) \\
The Mersenne primes $M_p$ are in one-to-one  &  The Euler primes $q$ are in one-to-one \\
correspondence with the even perfect numbers. & correspondence with the odd perfect numbers. \\
\hline
(E-2) The Mersenne primes $M_p$ satisfy & (O-2) The Euler primes $q$ satisfy \\
$M_p \equiv 3 \pmod 4$. & $q \equiv 1 \pmod 4$. \\
(Trivial) & (Trivial) \\
\hline
(E-3) The exponent $s = \nu_{M_p}(M)$ is one. & (O-3) The exponent $k = \nu_{q}(N)$ is one. \\
(Trivial) & (Sorli's Conjecture, 2003 \cite{Sorli}) \\
\hline
(E-4) If $M \equiv 0 \pmod 2$ is perfect, & (O-4) If $N \equiv 1 \pmod 2$ is perfect, \\
then given the Euclidean form & then given the Eulerian form \\
$M = {2^{p - 1}}(2^p - 1) = \displaystyle\prod_{i=1}^{2}{{p_i}^{\alpha_i}}$, & $N = {q^k}{n^2} = \displaystyle\prod_{j=1}^{\omega(N)}{{q_j}^{\beta_j}}$,\\
then ${{p_i}^{\alpha_i}\sigma({p_i}^{\alpha_i})}/M = i$, & then ${{q_j}^{\beta_j}\sigma({q_j}^{\beta_j})}/N \le 2/3 < j$, \\
for $i = 1, 2$. & for all $j$, $1 \le j \leq \omega(N)$. \\
(Observation, Dris 2011) & (Theorem, Dris 2008 \cite{Dris5}) \\
\hline
(E-5) There are infinitely & (O-5) There do not exist \\
many even perfect numbers. & any odd perfect numbers. \\
(EPN Conjecture) & (OPN Conjecture) \\
\hline
(E-6) The density of even & (O-6) The density of odd \\
perfect numbers is zero. &  perfect numbers is zero. \\
(Kanold 1954 \cite{Kanold}) & (Kanold 1954 \cite{Kanold}) \\
\hline
(E-7) $1 < I(M_p) \leq \frac{8}{7}$ for $p \geq 3$. & (O-7) $1 < I(q^k) < \frac{5}{4}$ for $q \geq 5$. \\
$\frac{7}{4} \leq \frac{2}{I(M_p)} = I(\frac{M_p + 1}{2}) < 2$ & $\frac{8}{5} < \frac{2}{I(q^k)} = I(n^2) < 2$ \\
In particular, & In particular, \\
$\frac{8}{7} < \sqrt{\frac{7}{4}} < I(\sqrt{\frac{M_p + 1}{2}}) < 2$. & $\frac{5}{4} < \sqrt{\frac{8}{5}} < I(n) < 2$.  \\
\hline
(E-8) $\frac{M_p + 1}{2} < M_p$ for $p \geq 3$. & (O-8) Conjecture: $q^k < n$ (Dris 2008 \cite{Dris5}). \\
\hline
(E-9) An even perfect number $M$& (O-8) An odd perfect number $N$ \\
has exactly two distinct prime factors. & has more than two distinct prime factors. \\
(i.e., $\omega(M) = 2$)  & (In fact, we know that  $\omega(N) \geq 9$ \cite{Nielsen2}.) \\
\hline
(E-10) $\gcd({2^p - 1}, 2^{p - 1}) = 1$ & (O-9) $\gcd(q^k, n^2) = \gcd(q, n) = 1$ \\
(Trivial) & (Euler) \\
\hline
\end{tabular}
\end{center}

\begin{rem}\label{Remark1}
We excluded $p_1 = 2$ from (E-7) because $\displaystyle\frac{M_{p_1} + 1}{2} = 2^{p_1 - 1} = 2^{2 - 1} = 2$ is squarefree.
\end{rem}

\begin{rem}\label{Remark2}
A proof of the inequality $q < n$ (which follows from $q^k < n$) has been announced in (Brown 2016) \cite{Brown}, (Starni 2017) \cite{Starni}, and (Dris 2017) \cite{Dris6}.  Essentially, Brown and Dris prove this claim by showing that Sorli's conjecture implies that $q^k < n$.
\end{rem}

In the next section, we give some \emph{known} relationships between the divisors of even and odd perfect numbers.  We will also discuss a conjectured relationship between the divisors $q^k$ and $n$ of an odd perfect number, which first appeared in the M.~Sc. thesis \cite{Dris5} that was completed in August of 2008.

\section{Inequalities Relating Divisors of Perfect Numbers}\label{Section2}
From Section \ref{Section1}, note that the heuristic (E-4), upon setting $Q = {2^p - 1}$, $K = 1$, and ${\bar{n}}^2 = 2^{p - 1}$, actually gives

$$\frac{\sigma(Q^K)}{{\bar{n}}^2} = \frac{\sigma(2^p - 1)}{2^{p - 1}} = \frac{(2^p - 1) + 1}{2^{p - 1}} = \frac{2^p}{2^{p - 1}} = 2,$$
and
$$\frac{\sigma({\bar{n}}^2)}{Q^K} = \frac{\sigma(2^{p - 1})}{2^p - 1} = \frac{2^{(p - 1) + 1} - 1}{(2 - 1)(2^p - 1)} = \frac{2^p - 1}{2^p - 1} = 1.$$

We state this result as our first lemma for this section.

\begin{lemma}\label{Lemma1}
If $M = {Q^K}{{\bar{n}}^2} = (2^p - 1)\cdot{2^{p - 1}}$ is an even perfect number given in Euclidean form, then we have the inequality
$$\frac{\sigma({\bar{n}}^2)}{Q^K} = 1 < 2 = \frac{\sigma(Q^K)}{{\bar{n}}^2}.$$
\end{lemma}

\begin{rem}\label{Remark3}
Except for the case of the first even perfect number $M = 6$ (as was pointed out in Remark \ref{Remark1}), the abundancy indices of the divisors of an even perfect number given in the Euclidean form $M = {Q^K}{{\bar{n}}^2} = (2^p - 1)\cdot{2^{p - 1}}$ (where the relabelling is done to mimic the appearance of the variables in the Eulerian form of an odd perfect number $N = {q^k}{n^2}$) satisfy the inequality
$$1 < I(Q^K) \leq \frac{8}{7} < \frac{7}{4} \leq I({\bar{n}}^2) < 2,$$
as detailed out in heuristic (E-7).  In particular, the inequality $I(Q^K) < I({\bar{n}}^2)$, together with Lemma \ref{Lemma1}, imply the inequality
$${\bar{n}}^2 < Q^K.$$
(In other words, we have the inequality
$$2^{p - 1} < 2^p - 1$$
where $p$ and $M_p = 2^p - 1$ are both primes.  Compare this with the inequality
$$q^k < n^2$$
for the divisors of an odd perfect number given in the Eulerian form $N = {q^k}{n^2}$ [see \cite{Dris3}, \cite{Dris5}].)
\end{rem}

The following result is taken from \cite{Dris3} and \cite{Dris5}.

\begin{lemma}\label{Lemma2}
If $N = {q^k}{n^2}$ is an odd perfect number given in Eulerian form, then we have the inequality
$$\frac{\sigma(q^k)}{n^2} \leq \frac{2}{3} < 3 \leq \frac{\sigma(n^2)}{q^k}.$$
\end{lemma}

Going back to the relabelling $Q = {2^p - 1}$, $K = 1$, and ${\bar{n}}^2 = 2^{p - 1}$ for an even perfect number $M = {Q^K}{{\bar{n}}^2} = (2^p - 1)\cdot{2^{p - 1}}$ given in Euclidean form, we now compute:

$$\frac{\sigma(Q^K)}{\bar{n}} = \frac{\sigma(2^p - 1)}{2^{\frac{p - 1}{2}}} = \frac{(2^p - 1) + 1}{2^{\frac{p - 1}{2}}} = \frac{2^p}{2^{\frac{p - 1}{2}}} = 2^{p - (\frac{p - 1}{2})} = 2^{\frac{p + 1}{2}},$$
and
$$\frac{\sigma(\bar{n})}{Q^K} = \frac{\sigma(2^{\frac{p - 1}{2}})}{2^p - 1} = \frac{2^{\frac{p - 1}{2} + 1} - 1}{(2 - 1)(2^p - 1)} = \frac{2^{\frac{p + 1}{2}} - 1}{2^p - 1}.$$

Now observe that
$$2^{\frac{p + 1}{2}} - 1 < 2^p - 1$$
since $1 < \frac{p + 1}{2} < p$, while we also have
$$2^{\frac{p + 1}{2}} \geq 4$$
because $p \geq 3$.  (Again, we excluded the first (even) perfect number $M = 6$ from this analysis because it is \emph{squarefree}.)

These preceding numerical inequalities imply that
$$\frac{\sigma(Q^K)}{\bar{n}} = 2^{\frac{p + 1}{2}} \geq 4 > 1 > \frac{2^{\frac{p + 1}{2}} - 1}{2^p - 1} = \frac{\sigma(\bar{n})}{Q^K}.$$

We state the immediately preceding result as our third lemma for this section.

\begin{lemma}\label{Lemma3}
If $M = {Q^K}{{\bar{n}}^2} = (2^p - 1)\cdot{2^{p - 1}}$ is an even perfect number given in Euclidean form (and $M \neq 6$), then we have the inequality
$$\frac{\sigma(Q^K)}{\bar{n}} \geq 4 > 1 > \frac{\sigma(\bar{n})}{Q^K}.$$
\end{lemma}

\begin{rem}\label{Remark4}
In particular, observe that the inequality
$$\frac{\sigma(\bar{n})}{Q^K} < 1$$
from Lemma \ref{Lemma3} implies that
$$\bar{n} < Q^K,$$
which, of course, \emph{trivially} follows from the inequality
$${\bar{n}}^2 < Q^K$$
in Remark \ref{Remark3}.

Likewise, compare the inequality (from Lemma \ref{Lemma3})
$$\frac{\sigma(\bar{n})}{Q^K} < \frac{\sigma(Q^K)}{\bar{n}},$$
for an even perfect number $M = {Q^K}{{\bar{n}}^2} = (2^p - 1)\cdot{2^{p - 1}}$ given in Euclidean form, with the inequality
$$\frac{\sigma(q^k)}{n} < \frac{\sigma(n)}{q^k},$$
for an odd perfect number $N = {q^k}{n^2}$ given in Eulerian form.  This second inequality was originally conjectured in \cite{Dris3-1}, and in fact, it has been recently shown (see \cite{Dris6}) to be equivalent to the conjecture $q^k < n$, which originally appeared in the M.~Sc. thesis \cite{Dris5}.

Details for the proof of the biconditional
$$q^k < n \Longleftrightarrow \sigma(q^k) < \sigma(n) \Longleftrightarrow \frac{\sigma(q^k)}{n} < \frac{\sigma(n)}{q^k}$$
are clarified in (Dris 2017 \cite{Dris6}).
\end{rem}

The next section will explain our motivation for pursuing a proof for the following conjecture:

\begin{conj}\label{Conjecture1}
If $N = {q^k}{n^2}$ is an odd perfect number given in Eulerian form, then the conjunction
$$\{k = \nu_{q}(N) = 1\} \land \{q^k < n\}$$
always holds.
\end{conj}

\section{On the Conjectures of Sorli and Dris Regarding Odd Perfect Numbers}\label{Section3}
We begin this section with a recapitulation of the two main conjectures on odd perfect numbers that have been mentioned earlier in this article.

\begin{conj}\label{Conjecture2}
The Descartes-Frenicle-Sorli conjecture predicts that if $N = {q^k}{n^2}$ is an odd perfect number given in Eulerian form, then
$$k = \nu_{q}(N) = 1$$
always holds.
\end{conj}

\begin{rem}\label{Remark5}
Dris gave a sufficient condition for Sorli's conjecture in \cite{Dris3}.  Some errors, however, were found in the initial published version of that article, and Dris had to retract his claim that the biconditional
$$k = \nu_{q}(N) = 1 \Longleftrightarrow n < q$$
always holds.  (The current published version of \cite{Dris3} contains a proof only for the one-sided implication $$n < q \Longrightarrow k = \nu_{q}(N) = 1.$$  
In two [continually] evolving papers [see \cite{Dris} and \cite{Dris2}], work is in progress to try to disprove the converse
$$k = \nu_{q}(N) = 1 \Longrightarrow n < q$$
and thereby get a proof for Conjecture \ref{Conjecture1}.)

Moreover, Acquaah and Konyagin \cite{AcquaahKonyagin} \emph{almost} disproves $n < q$ by obtaining the estimate $q < n\sqrt{3}$ under the assumption $k = \nu_{q}(N) = 1$.  (Since the contrapositive of the implication $n < q \Longrightarrow k = 1$ is $k > 1 \Longrightarrow q < n$, we know that Acquaah and Konyagin's estimate for the Euler prime $q$ implies that the inequality
$$q < n\sqrt{3}$$
holds unconditionally.)

With the recent proofs for the inequality $q < n$ by Brown \cite{Brown}, Starni \cite{Starni}, and Dris \cite{Dris6}, it remains to show that the Descartes-Frenicle-Sorli conjecture holds so as to prove that the Dris conjecture that $q^k < n$ for odd perfect numbers $q^k n^2$ with Euler prime $q$ is indeed true.

Curiously enough, the two papers \cite{Dris3} and \cite{Dris-2} by Dris are cited in OEIS sequence A228059 \cite{OEIS-A228059}, whose description is reproduced below:

Odd numbers of the form $r^{1+4l}s^2$, where $r$ is prime of the form $1+4m$, $s > 1$, and $\gcd(r,s) = 1$ that are \emph{closer to being perfect than previous terms}.

\emph{Coincidentally}, the ``Euler prime'' of the first $9$ terms in this OEIS sequence all have exponent $1$:

$$45 = {5}\cdot{3^2}$$
$$405 = {5}\cdot{3^4}$$
$$2205 = {5}\cdot{(3\cdot7)}^2$$
$$26325 = {13}\cdot{({3^2}\cdot{5})}^2$$
$$236925 = {13}\cdot{({3^3}\cdot{5})}^2$$
$$1380825 = {17}\cdot{({3}\cdot{5}\cdot{19})}^2$$
$$1660725 = {61}\cdot{({3}\cdot{5}\cdot{11})}^2$$
$$35698725 = {61}\cdot{({3^2}\cdot{5}\cdot{17})}^2$$
$$3138290325 = {53}\cdot{({3^4}\cdot{5}\cdot{19})}^2.$$
\end{rem}

\begin{conj}\label{Conjecture3}
Dris's conjecture states that if $N = {q^k}{n^2}$ is an odd perfect number given in Eulerian form, then
$$q^k < n$$
always holds.
\end{conj}

\begin{rem}\label{Remark6}
Prior to the paper \cite{AcquaahKonyagin} by Acquaah and Konyagin, and the data from OEIS sequence A228059 \cite{OEIS-A228059} as detailed out in Remark \ref{Remark5}, the only heuristic available to justify Dris's conjecture that $q^k < n$ was the inequality
$$I(q^k) < \frac{5}{4} < \sqrt[3]{2} < \sqrt{\frac{8}{5}} < I(n).$$
(See the paper \cite{Dris3} for a proof.)  In particular, the \emph{heuristic justification} is that the divisibility constraint $\gcd(q^k, n) = \gcd(q, n) = 1$ appears to induce an 
``ordering property'' between $q^k$ and $n$ via an appropriate inequality between their abundancy indices.  That is, Dris expects his conjecture $q^k < n$ to follow from the last inequality above, in the sense that the inequality $q^k < n^2$ appears to have followed from the related inequality
$$I(q^k) < \frac{5}{4} < \sqrt{2} < \frac{8}{5} < I(n^2).$$

Additionally, note that all of the $8$ terms (apart from the first one) in the OEIS sequence mentioned in Remark \ref{Remark4} satisfy Dris's conjecture. 
\end{rem}

\section{Conclusion}
To conclude, a recent e-mail correspondence of the author with Brian D. Beasley of Presbyterian College revealed the following information, quoted verbatim from page 25 of \cite{Beasley}:

``Before proceeding with Euler's proof, we pause to note that his result was not quite what Descartes and Frenicle had conjectured, as they believed that $k=1$, but it came very close. In fact, current research continues in an effort to prove $k=1$. For example, Dris has made progress in this direction, although his paper refers to Descartes' and Frenicle's claim (that $k=1$) as Sorli's conjecture; Dickson has documented Descartes's conjecture as occurring in a letter to Marin Mersenne in 1638, with Frenicle's subsequent observation occurring in 1657.''

It might be wise (at this point) to delve deeper into this little bit of history in mathematics, to attempt to answer the particular question of whether Descartes and Frenicle used \emph{similar} or totally different methods to arrive at what we have come to call as Sorli's conjecture on odd perfect numbers.  Perhaps they \emph{both} used methods similar to the ones used in this article - who knows?  Besides, Mersenne's predictions for succeeding primes $p$ for which $2^p - 1$ turned out to be a Mersenne prime were already stunning as they were.  Did Mersenne use an \emph{algorithm}, for testing primality of Mersenne prime-number candidates, that remains unknown to the rest of us to this day?  Only time can tell.

% For alignments use AmS-LaTeX constructions not \eqnarray.

%% - theorems and proofs
%%\begin{thm}[optional text]
% The optional material will be typeset as part of the theorem heading
%%\end{thm}

%%\begin{proof}[Optional proof heading]
% the proof
%%\end{proof}
% An end-of-proof symbol (open box) will be typeset at the
% end of the proof.

\ack % or \acks
% Put acknowledgements here
The author would like to thank Professor Carl Pomerance for pointing out the relevance of the paper \cite{AcquaahKonyagin}, as well as Peter Acquaah for helpful e-mail exchanges on the topic of odd perfect numbers.  Lastly, he would like to express his gratitude to his research collaborator at Far Eastern University, Keneth Adrian P. Dagal.

% alteratively, bibliographies prepared with BibTeX can be included by
% means of the following commands
%\bibliographystyle{srtnumbered}
%\bibliography{mybib}

\begin{thebibliography}{20}
\bibitem{AcquaahKonyagin}
P.~Acquaah, S.~Konyagin, On prime factors of odd perfect numbers, {\it Int. J. Number Theory}, {\bf 08} (2012), 1537, doi:\url{https://dx.doi.org/10.1142/S1793042112500935}.
\bibitem{Beasley}
B.~D.~Beasley, Euler and the ongoing search for odd perfect numbers, ACMS 19th Biennial Conference Proceedings, Bethel University, May 29 to Jun. 1, 2013, \url{http://godandmath.files.wordpress.com/2013/07/acms-2013-proceedings.pdf}, pages 21-31.
\bibitem{BroughanDelbourgoZhou}
K.~A.~Broughan, D.~Delbourgo and Q.~Zhou, Improving the Chen and Chen result for odd perfect numbers, {\it Integers}, {\bf 13} (2013), Article $\#$A39, \url{https://www.emis.de/journals/INTEGERS/papers/n39/n39.pdf}, ISSN 1867-0652.
\bibitem{Brown}
P.~A.~Brown, A partial proof of a conjecture of Dris, preprint, \url{https://arxiv.org/abs/1602.01591}.
\bibitem{ChenChen}
F.~J.~Chen, Y.~G.~Chen, On the index of an odd perfect number, {\it Colloq. Math.}, {\bf 136} (2014), 41-49, doi:\url{https://dx.doi.org/10.4064/cm136-1-4}.
\bibitem{Dris}
J.~A.~B.~Dris, New results for the Descartes-Frenicle-Sorli conjecture on odd perfect numbers, {\it Journal for Algebra and Number Theory Academia}, {\bf 6} (2016), issue 3, 95-114.
\bibitem{Dris2}
J.~A.~B.~Dris, New results for Sorli's conjecture on odd perfect numbers - Part II, (Aug. 2013), preprint, \url{https://arxiv.org/abs/1303.2329}.
\bibitem{Dris-2}
J.~A.~B.~Dris, A short "proof" for Sorli's conjecture on odd perfect numbers, (Oct. 2014), preprint, \url{https://arxiv.org/abs/1308.2156}.
\bibitem{Dris3}
J.~A.~B.~Dris, The abundancy index of divisors of odd perfect numbers, {\it J. Integer Seq.}, {\bf 15} (Sep. 2012), Article 12.4.4, \url{https://cs.uwaterloo.ca/journals/JIS/VOL15/Dris/dris8.html}, ISSN 1530-7638.
\bibitem{Dris3-1}
J.~A.~B.~Dris, Some new results on bounds for the abundancy indices of the components of odd perfect numbers, (Mar. 2011), preprint, \url{https://arxiv.org/pdf/1103.1090v1.pdf}.
\bibitem{Dris4}
J.~A.~B.~Dris, Perfect numbers - On Mersenne and Euler primes, \url{https://mathoverflow.net/questions/48953}, asked on 12/10/2010.
\bibitem{Dris5}
J.~A.~B.~Dris, Solving the Odd Perfect Number Problem: Some Old and New Approaches, M.~Sc. thesis, De La Salle University, Manila, Philippines, 2008, \url{https://arxiv.org/abs/1204.1450}.
\bibitem{Dris6}
J.~A.~B.~Dris, On a curious biconditional involving the divisors of odd perfect numbers, preprint, \url{https://arxiv.org/pdf/1309.0906v19.pdf}.
\bibitem{Kanold}
H.-J.~Kanold, Uber die Dichten der Mengen der vollkommenen und der befreundeten Zahlen, {\it Math. Z.}, {\bf 61} (1954), 180-185, doi:\url{https://dx.doi.org/10.1007/BF01181341}.
\bibitem{Nielsen1}
P.~P.~Nielsen, Odd perfect numbers, Diophantine equations, and upper bounds, {\it Math. Comp.}, {\bf 84} (2015), 2549-2567, doi:\url{https://doi.org/10.1090/S0025-5718-2015-02941-X}.
\bibitem{Nielsen2}
P.~P.~Nielsen, Odd perfect numbers have at least nine distinct prime factors, {\it Math. Comp.}, {\bf 76} (2007), 2109-2126, doi:\url{https://dx.doi.org/10.1090/S0025-5718-07-01990-4}.
\bibitem{OchemRao}
P.~Ochem, M.~Rao, Odd perfect numbers are greater than ${10}^{1500}$, {\it Math. Comp.}, {\bf 81} (2012), 1869-1877, doi:\url{https://dx.doi.org/10.1090/S0025-5718-2012-02563-4}.
\bibitem{OEIS-A228059}
T.~D.~Noe, Online Encyclopedia of Integer Sequences, sequence A228059, \url{http://oeis.org/A228059}, last accessed: 10/12/2013.
\bibitem{Starni}
P.~Starni, On Dris conjecture about odd perfect numbers, preprint, \url{https://arxiv.org/abs/1706.02144}.
\bibitem{Sorli}
R.~M.~Sorli, Algorithms in the Study of Multiperfect and Odd Perfect Numbers, Ph.~D. Thesis, University of Technology, Sydney, 2003, \url{https://epress.lib.uts.edu.au/research/handle/10453/20034}.
\end{thebibliography}

\end{document}